\documentstyle[amscd,amssymb,12pt]{amsart}
\input xypic \xyoption{curve}
\input epsf

\def\hcorrection#1{\advance\hoffset by #1 }
\def\vcorrection#1{\advance\voffset by #1 }

\vcorrection{-.7in}
\hcorrection{-.5in}

\newcommand{\B}[1]{{\bold#1}} 
\newcommand{\C}[1]{{\cal#1}} 

\theoremstyle{plain}
\newtheorem{th}{Theorem}[section]

\newtheorem{lem}{Lemma}[section]
\newtheorem{prop}{Proposition}[section]

\theoremstyle{definition}
\newtheorem{defin}{Definition}[section]

\theoremstyle{definition}

\theoremstyle{remark}
\newtheorem{rem}{Remark}[section]

\numberwithin{equation}{section}

\begin{document}

\pagestyle{plain}
\addtolength{\footskip}{.3in}

\title[Coefficients in deformation quantization]
{A combinatorial approach to coefficients in deformation quantization}
\author{Lucian M. Ionescu}
\address{Department of Mathematics, Illinois State University, IL 61790-4520
}
\email{lmiones@@ilstu.edu}
\date{04/20/2004}

\begin{abstract}
Graph cocycles for star-products are investigated from 
the combinatorial point of view, 
using Connes-Kreimer renormalization techniques.

The Hochschild complex, 
controlling the deformation theory of associative algebras,
is the ``Kontsevich representation'' of a DGLA of graphs 
coming from a pre-Lie algebra structure defined by graph insertions
(Gerstenhaber composition with Leibniz rule).

Properties of the dual of its UEA (an odd parity analog of Connes-Kreimer Hopf algebra),
are investigated in order to find solutions of the deformation equation.
The solution of the initial value deformation problem, at tree-level, is unique.
For linear coefficients the resulting formulas are relevant to the Hausdorff series.
\end{abstract}

\maketitle
\tableofcontents

\clearpage


\section{Introduction}\label{S:intro}
Kontsevich solution for deformation quantization of Poisson manifolds \cite{Kon1}
is based on a specific QFT \cite{CF},
yet universal: use one quantization to find them all.
The coefficients of the star-product,
expressed as integrals on some configuration spaces,
are the corresponding amplitudes.

In \cite{pqft} it was shown that the coefficients represent a cocycle
of the cobar construction of the dg-coalgebra of Kontsevich graphs,
and it was conjectured that there should exist a direct combinatorial approach,
as an alternative to the above analytical approach.

In this article the deformation problem for associative algebras is
pulled back at the level of graphs.
As a generalization,
generic associative deformations of associative algebras in a generic direction,
not necessarily Poisson,
are also considered.
The corresponding cocycles are ``unitary characters'',
i.e. satisfying $W_c^{-1}=\bar{W}_c$ (inverse relative to convolution).
The ``initial value (deformation) problem'' at tree-level (no loops)
has a unique solution.
In the case of linear coefficients,
e.g. for $g^*$, with $g$ a Lie algebra,
the solution is expected to correspond to the Hausdorff series,
which is the ``forefather'' of all quantizations.

This is possible since the Kontsevich rule (a kind of a dual Feynman rule),
is a cohomological pairing,
allowing to pullback Gerstenhaber composition from the Hochschild complex.

%
Besides its value as a self contained approach (at the mathematical level)
to the quantization problem,
the present approach is relevant to the understanding of renormalization techniques.
Not surprisingly, since the coefficients may be interpreted as Feynman integrals,
the ``correction analysis'' of the star-product terms in order to yield
associativity is similar to BPHZ renormalization via ``counter-terms method'' \cite{Kre1,Kre2,CKren1}.
Connections between deformation theory and renormalization are expected to emerge.

To make the article accessible to a wider audience,
the next section is a brief recall on deformation quantization
and the associated combinatorial problem
(Students could perhaps start with sections \ref{S:dq}, \ref{S:ca} and \ref{S:scr}).

Section \ref{S:gc} continues with a short recall on graph cohomology from \cite{pqft},
the reduction from formality to star-products leading to 
section \ref{S:sp} where the relevant graph operations are discussed.
Examples are computed in section \ref{S:coa} (see also Appendix),
which rephrases the ``correction analysis'' of \cite{Kath} 
in the context of the Hopf algebra of graphs.

The reader in deformation theory and renormalization 
is invited to skim through section \ref{S:sp} before considering 
where the correspondence between associativity and the unitarity condition
for the coefficients is explained,
before considering the main point stated in \ref{S:dgla}.

The relation with renormalization is briefly considered 
as part of the concluding section.

\vspace{.1in}
{\bf Acknowledgments}
Comments are most welcome!

\section{Background in deformation quantization}\label{S:dq}
In classical physics observables like position $x$ and momentum $p$ commute. 
In the quantum world they don't. 
They are usually modeled by linear operators on Hilbert spaces, 
satisfying the celebrated Heisenberg uncertainty relation: $x p - p x=hI$. 

The ``Quantization Problem'' consists in finding a correspondence 
between classical and quantum variables, 
in order to pass (and benefit) from a classical description 
to the quantum description. 
A broad and ``simple'' approach, 
the so called deformation quantization approach, 
is based on the idea of keeping the observables ``as is'' 
while deforming the commutative product into 
a non-commutative one: a ``star-product'' satisfying: 
$x p - p x=hI$.

Recall that the state space of a classical system is modeled by a 
Poisson manifold $(X,\alpha)$, 
where $\alpha$ denotes a two tensor: the Poisson structure. 
Assume for simplicity that $X=R^d$. 
In the case of a constant coefficients Poisson bracket, 
there is a nice exponential (Moyal's) formula for such a star-product: 
\begin{equation}\label{E:Moyal}
f \star g = exp(h\alpha)(f\otimes g),
\end{equation}
where $f,g$ denote functions on $X$.

In 1997, in a ground breaking paper \cite{Kon1}, 
Kontsevich devised a formula for the star-product in the general case (variable coefficients): 
$$f\star g=\sum_n\frac{h^n}{n!}\sum_{\Gamma\in G_n}W_\Gamma U_{\Gamma}(\alpha^n), \qquad
W_\Gamma=\frac1{n!(2\pi)^{2n}}\int_{H_n}\prod_{k=1}^n(d\phi_{e_k^1}\wedge d\phi_{e_k^2}^h).$$
The above star-product formula is to be thought of as a Taylor series,
and where the terms, 
due to the non-commutativity of the bidifferential operators $U(\alpha)$ involved, 
correspond to a certain class of graphs $G_n$ (not just uni-valent trees). 
The correspondence between graphs and bidifferential operators can be thought 
as an instance of a ``Feynman rule'' from quantum field theory (QFT).
It allows the use of such a graphical representation,
as a tool to tame the combinatorics lurking in the background.

The difficult part in the Kontsevich formula for the star-product 
is to come up with a formula for the coefficients $W$,
so that the product would be associative. 
The analytical solution of Kontsevich for the coefficients $W$ 
is based on integrals over some configuration spaces $H_n$.
Their interpretation as Feynman integrals was given in \cite{CF}.
In this sense Kontsevich solution uses {\em one} quantization
(with underlying Feynman path integral method)
to provide the adequate weights for the quantization
of {\em all} Poisson manifolds.

\vspace{.1in}
\hspace{1.8in}\epsfbox{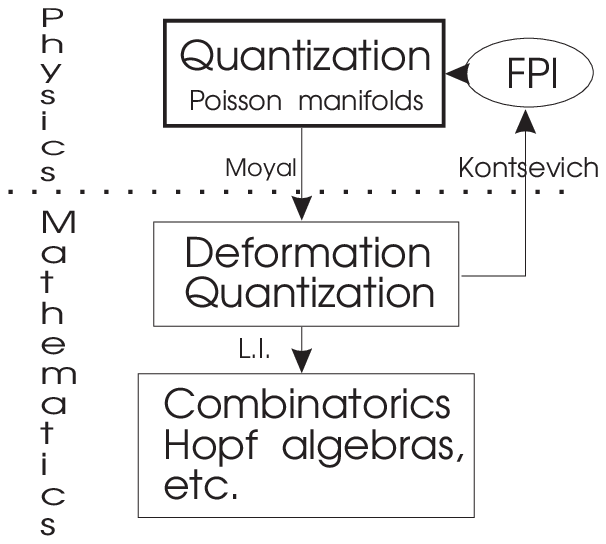}

\vspace{.1in}
What about a ``pure'' mathematical solution?
A tentative combinatorial approach is the primary concern of 
the present author's work.

The condition corresponding to associativity can be rewritten 
as an integral over the boundary of the compactification 
of the above mentioned spaces. 
It was noticed in \cite{pqft} that the various portions of 
the boundary in this multi-facet boundary correspond 
to a certain ``Forest Formula'' used in renormalization,
in order to tame the infinities appearing when computing 
the physicist's version of the Feynman integrals. 
With the same purpose, 
Kreimer \cite{Kre1} introduced on the vector space of graphs 
the structure of a coproduct (see also \cite{CKren1}). 
If a product ``merges'' the factors to get a result, the coproduct:
$$\Delta\Gamma=\sum_{\gamma\subset\Gamma\to \gamma'}\gamma\otimes\gamma',$$
encodes the internal structure of the graph $\Gamma$, 
representing a list of pairs consisting of a subgraph together 
with the result of collapsing it in the original bigger graph. 
It may be thought of as a list of ways to build the bigger graph from two pieces. 
This simple technique corresponding to an almost obvious idea (well known in combinatorics), 
clarified the renormalization problem (roughly speaking). 
The present author used this technique in connection with the formality morphism \cite{pqft}, 
and obtained that the coefficients are (co)cycles determined by 
a closely related coproduct (see also \cite{cfg}).
 
This opened a new research avenue of finding a combinatorial formula 
for these (or another family of) coefficients,
satisfying the cocycle condition: $\delta W=0$.
It would not only simplify the computations of 
the deformation quantization formula, 
but would have applications in the theory of Lie algebras. 
The Campbell-Baker-Hausdorff formula, 
which may be thought of as the forefather of all quantizations, 
is in need of a ``reasonable explicit expression'' for its coefficients (\cite{Kath}, p.7). 
Its terms represent various ways of bracketing Lie algebra elements,
corresponding to $U_\Gamma(\alpha^n)$,
where $\Gamma$ is a binary tree with $n$ vertices,
and $\alpha$ is the Lie bracket.

\section{Graph cocycles}\label{S:gc}
Consider the coaugmented counital dg-coalgebra of 
Kontsevich graphs $(H,\cdot, u, \Delta, \eta, d)$,
with $H=k\C{G}$ the k-vector space generated by
orientation classes of labeled directed graphs $[\Gamma,l]\in\C{G}$.
For additional details (not essential for what follows),
we direct the reader to \cite{pqft,cfg}.
The subset $G_{n,m}^l$ of $\C{G}$ will consist in graphs with
$n$ internal vertices, $m$ boundary vertices and $2n+m-2+l$ edges,
with no loops (including no circuits).

\hspace{2in}\epsfbox{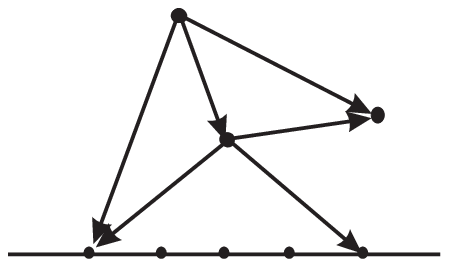}

\hspace{2.3in} $n=3,\ m=5,\ l=-3$

\vspace{.2in}
Recall that the cocycles condition refers to the dual $\delta=D^*$ 
of the cobar total differential:
$$D=D_d+D_{\Delta_b},$$
where $D_d$ and $D_{\Delta_b}$ are coderivations on the reduced tensor coalgebra $T(\bar{H})$,
extending the graph differential:
\begin{equation}\label{E:fgdiff}
d([\Gamma,l])=\sum_{e\in \Gamma^{(1)}_{int}} [\Gamma/e,l/e], 
\quad |\Gamma^{(1)}_{int}|>0,
\end{equation}
and the reduced coproduct $\bar{\Delta}=\Delta_b$:
\begin{equation}\label{E:fgcoprod}
\Delta_b[\Gamma,l]=\sum_{\gamma\to \Gamma\to \gamma', \ \gamma\cap\partial\Gamma\ne \emptyset}
[\gamma,l|\gamma]\otimes[\gamma',l/\gamma].
\end{equation}
The above sum is over all {\em normal subgraphs} of $\Gamma$ meeting the boundary $\partial\Gamma$,
i.e. such that collapsing $\gamma$ to a (boundary) vertex yields a graph from the given class $\C{G}$.
%

\hspace{.5in} \epsfbox{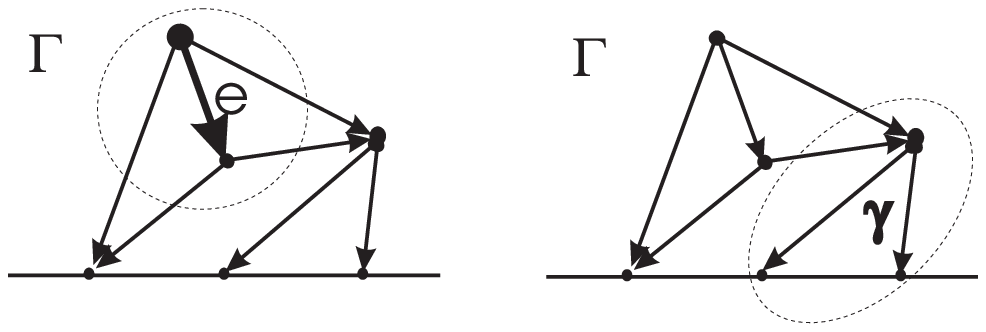}

We begin by introducing a special subclass of graphs.
\begin{defin}
Let $b_n^L\in G_{n,2}$ denote the orientation class
of the labeled {\em left Bernoulli graphs} (see also \cite{Kath}):
$$\diagram
1\dto_{l_1}\rto^{r_1} & 2\dlto_{l_2}\qquad ... & n\dllto_{l_n}\dto^{r_n}\\
L & & R.
\enddiagram$$
Correspondingly, the {\em right Bernoulli graph} $b_n^R$ 
is built out of wedges with all right legs grounded to the boundary point $R$.
\end{defin}
It may be obtained from the binary tree corresponding to the right parenthesized product of
$n$ letters,
by collapsing all left descending terminals to the boundary point $L$.
The unique terminal which is a right descendent will be denoted by $R$.

A few properties of Bernoulli graphs are listed next.
\begin{lem}
Let $b_n$ denote a left (right) Bernoulli graph. Then
(i) $b_n$ are primitive elements of the Hopf algebra of graphs.
(ii) $Db_n=db_n$.
(iii) The graph differential of a Bernoulli graph is:
$$n\ even:\ db_n=0, \qquad n\ odd:\ db_n=-b_{n-1}.$$
\end{lem}
\begin{pf}
The graphs $b_n$ are ``simple'', i.e. their only normal subgraphs are the trivial ones.
As a consequence $\bar{\Delta}(b_n)=0$, 
and (i) is established.
Rephrasing (i) gives (ii).

In order to prove (iii), 
note that $n-2$ equals the number of internal edges,
and that by collapsing an internal edge $e$, 
essentially the ``previous'' Bernoulli graphs is obtained:
$b_n/r_k=(-1)^k b_{n-1}$.
The claim follows after taking into account the orientation
which determines the alternating sign
(by definition $db_2=0$).
\end{pf}
%

Graph cocycles $W$ of the cobar construction applied to the above dg-coalgebra
determine L-infinity morphisms $U$ \cite{cfg}.
If $U$ is such a morphism and $\alpha$ is a Poisson structure, 
thus associative and satisfying Maurer-Cartan equation $d\alpha+[\alpha,\alpha]=0$, with $d=0$,
then:
$$\tilde{\alpha}=U(exp^\alpha)=\sum_{n=0}^\infty \frac1{n!}U(\alpha\wedge...\wedge\alpha)$$
satisfies the corresponding Maurer-Cartan equation (see \cite{Kon1}, p.12).
The corresponding restriction to a product on $A$ (Poisson algebra of functions):
$$f\star g=U(e^\alpha)(f\otimes g)$$
is associative, and thus yields a star-product.

The above ``philosophy'' can be applied to the Campbell-Baker-Hausdorff formula,
viewed as a deformation quantization of the abelian structure underlying the vector space structure,
in the direction of the Lie bracket (so called CBH-quantization \cite{Kath}):
$$H(x,y)=x+y+\epsilon \frac12[x,y]+\epsilon^2 \frac1{12}([x,[x,y]]+[[x,y],y])+....$$
The graphs involved are the {\em sym-admissible graphs} \cite{Kath},
for which the corresponding prime graphs are the {\em Lie admissible graphs} $\C{L}$.

It turns out that since $ad:g\to Der(g)$ is a Lie algebra morphism,
then the star-product is determined by the restriction to $x^n\star y$ (see \cite{Kath}),
for which the corresponding graphs are the Bernoulli graphs.

In what follows we will study the associativity condition from the point of view of
the above cocycle condition.

\section{Star-products}\label{S:sp}
Given $\C{U}$ and a system of weights $W$,
define the following product:
\begin{equation}\label{E:star-product}
f\star g=\sum_{n\ge 0}\epsilon^n\sum_{\Gamma\in G^0_{n,2}}W_\Gamma U_\Gamma(\alpha^n)(f,g),
\quad f,g\in A,
\end{equation}
where $A$ is the associative algebra $C^\infty(R^d)$ for simplicity.

From this point on, 
we will consider directed graphs with a fixed order on the boundary vertices,
and without additional labels, 
since the same polyvector field $\alpha$ will be attached to all its internal vertices,
and the labeling of the edges, used to define the ``state-sum'',
is here irrelevant.
The missing symmetry factors represent the price paid for this simplification.

The graphs with $n=0,1,2$ internal vertices from $G_{n,2}^0$ are 
(``prime'') Bernoulli graphs $b_0, b_1, b_2^{L/R}$,
or ``products'' of Bernoulli graphs $b_1^2$ (see \ref{D:prod}).

\vspace{.2in}
\hspace{1.5in} \epsfbox{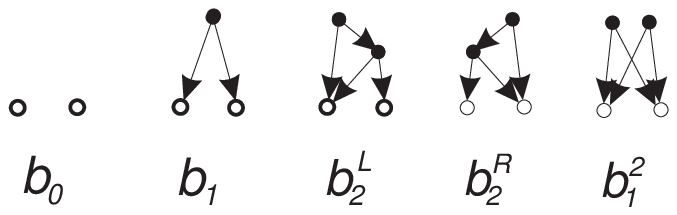} 

The product alluded to above, introduced by \cite{Kath}, p.22, 
corresponds to the disjoint union of graphs 
with the same boundary/external structure.
\begin{defin}\label{D:prod}
If $\Gamma_1\in G_{n,m},\Gamma_2\in G_{n',m}$ define
$\Gamma_1\Gamma_2\in G_{n+n',m}$ as the graph
obtained from the disjoint union of the two graphs
by identification of the boundary vertices.
Extend the product to all pairs of graphs,
so that all other products are zero.
\end{defin}
For example $b_1^2=b_1 b_1$, therefore justifying the notation used.
Note that graphs have a unique factorization relative to the above product.
The graph is called {\em prime} if it is not a non-trivial product of graphs.

\vspace{.1in}
A few graphs from $G_{n,3}$ will be used in examples.
The graphs with $n=1$ internal vertices are denoted by $b_1^{L/M/R}$:

\hspace{1.6in} \epsfbox{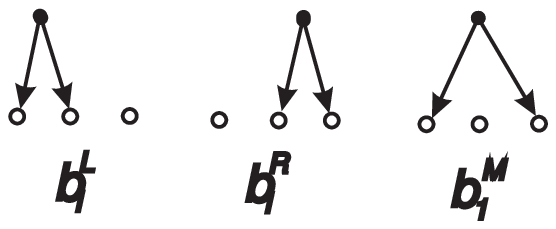}

The following graphs for $n=2$ will be denoted by
$t_2^{L/R}$, $c_2$ (prime graphs) and $c_2^{L/R}$.

\hspace{1in} \epsfbox{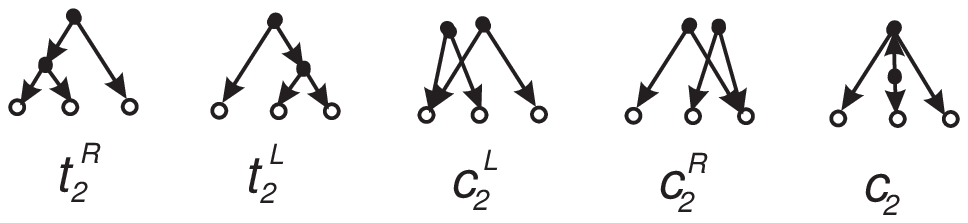}

The formula \ref{E:star-product} produces a star-product, 
i.e. it is associative, 
iff for any $n\ge 0$ and functions $f,g,h$:
$$\sum_{k+l=n}\sum_{\Gamma_1\in G_{k,2}^0}\sum_{\Gamma_2\in G_{l,2}^0}
W_{\Gamma_1}W_{\Gamma_2}[\C{U}_{\Gamma_1}(\C{U}_{\Gamma_2}(f,g),h)-
\C{U}_{\Gamma_1}(f,\C{U}_{\Gamma_2}(g,h))]=0.$$
Using Gerstenhaber composition \cite{GS}, 
the above equations can be more compactly written as:
\begin{equation}\label{E:assoc1}
0=\star\circ\star=\sum_{n\ge 0}\epsilon^n
[\sum_{k+l=n}\sum_{\Gamma_1\in G_{k,2}^0}\sum_{\Gamma_2\in G_{l,2}^0}
W_{\Gamma_1}W_{\Gamma_2}
(\C{U}_{\Gamma_1}\circ \C{U}_{\Gamma_2})].
\end{equation}
The above composition of operators can be replaced with the operator
associated to a certain product of graphs.
\begin{defin}\label{D:star1}
The {\em composition of graphs} is defined as follows:
$$\Gamma_1\circ\Gamma_2=\sum_{i=1}^m (-1)^{(i-1)(m'-1)}\Gamma_1\circ_i\Gamma_2,
\quad \Gamma_1\in G_{n,m}^l, \ \Gamma_2\in G_{n',m'}^{l'},$$
where $\circ_i$ is the insertion of $\Gamma_2$ at the $i^{th}$
boundary vertex of $\Gamma_1$ using ``Leibniz rule''
i.e. summing over all possible graphs where the ``$i^{th}$ legs'' of $\Gamma_1$ lend on 
vertices of $\Gamma_2$.
\end{defin}
If $[I_i]$ denotes the set of incoming edges at the $i^{th}$ boundary point
of $\Gamma_1$ and $[n']$ denotes the set of internal vertices of $\Gamma_2$,
then the ``Leibniz rule'' at the $i^{th}$ vertex yields:
$$\Gamma_1\circ_i\Gamma_2=\sum_{f:[I_i]\to [n']}\Gamma_1\circ_i^f\Gamma_2,$$
where $\circ_i^f$ denotes the operation of 
replacing the vertex $i$ by the set $[n']$ and 
the edges $e=(v\to i)\in [I_i]$ by $(v\to f(e))$.
The insertion operation is compatible with the product of graphs.
\begin{lem}
The linear extension of the insertion operation $\circ_i$ is an algebra morphism:
$$(\Gamma_1\Gamma_2)\circ_i\Gamma=(\Gamma_1\circ_i\Gamma)(\Gamma_2\circ_i\Gamma).$$
\end{lem}
\begin{pf}
The lemma amounts to the fact that gluing data $f:[I^1_i]\cup[I^2_i]\to [n']$ is 
determined by $f_j:[I^j_i]\to [n'], j=1,2$.
\end{pf}
For the case $m=2$ (see \cite{Kath}, p.16):
$$\Gamma_1\circ \Gamma_2=\Gamma_1\circ_L\Gamma_2-\Gamma_2\circ_R\Gamma_1,$$
where the indices $1$ and $2$ where replaced by the more suggestive ones $L$ and $R$,
representing the left and right boundary points respectively.
Note that in this case the product of two graphs from $G_{k,2}^0$ and $G_{l,2}^0$
is a graph of $G_{k+l,3}^{-1}$.

A few examples are in order:
\begin{equation}\label{E:comp1}
b_0\circ b_0=0, \quad [b_0,b_1]=0,\quad 
b_1\circ b_1=t_2^R-t_2^L+c_2^L-c_2^R.
\end{equation}
The last computation of $b_1\circ b_1$ is represented below:

\vspace{.1in}
\hspace{.01in}\epsfbox{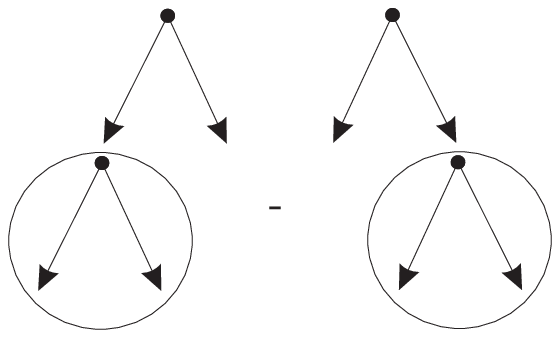}
\hspace{.2in}\epsfbox{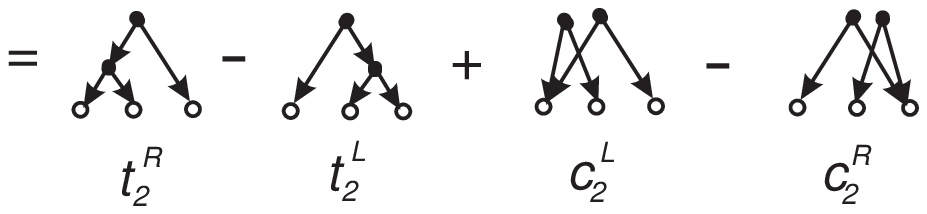}

%
This product is essentially the pullback of Gerstenhaber composition via the
``Feynman rule'' $\C{U}$.
It is similar to the pre-Lie operation on graphs defined in \cite{Kre2}, p.27.,
except insertions are allowed only at boundary vertices.

The {\em associator} of graph composition is:
$${\bf \alpha}(\Gamma_1,\Gamma_2,\Gamma_3)=(\Gamma_1\circ\Gamma_2)\circ\Gamma_3-
\Gamma_1\circ(\Gamma_2\circ\Gamma_3).$$
\begin{defin}\label{D:assoc}
The associated graded (Lie) bracket is defined by:
\begin{equation}\label{E:Liebracket}
[\Gamma_1,\Gamma_2]=\Gamma_1\circ\Gamma_2-(-1)^{(m-1)(m'-1)}\Gamma_2\circ\Gamma_1,
\quad \Gamma_1\in G_{n,m}, \ \Gamma_2\in G_{n',m'}.
\end{equation}
\end{defin}
%
\begin{lem}
The associator satisfies the following equation:
$${\bf \alpha}(\Gamma_1,\Gamma_2,\Gamma_3)=-\alpha(\Gamma_1,\Gamma_3,\Gamma_2),$$
therefore the composition of graphs is a pre-Lie composition
with associated Lie bracket $[\ ,\ ]$.
\end{lem}
\begin{pf}
It is well know that the above antisymmetry of the associator
is equivalent with the associated bracket satisfying 
the Jacobi identity \cite{HDGLA}, p.5.

A proof similar to the one for the Gerstenhaber composition
of Hochschild cochains $f\circ_i g$ establishes the above equation.
The only difference is that here $\Gamma_1\circ_i\Gamma_2$,
which is the analog of $f\circ_i g$, is a sum of graphs.
Otherwise, the signs and symmetry properties are identical.
\end{pf}
\begin{lem}\label{L:loops}
The composition of two graphs without circuits is a sum of graphs without circuits.
\end{lem}
\begin{pf}
Since all the arrows connecting the first graph enter the second graph,
$\Gamma_1\circ\Gamma_2$ has no circuits.
\end{pf}
Consider the basis in $H^*$ dual to $G$,
consisting in functionals $\delta_\Gamma, \Gamma\in G$.
The composition of graphs induces a ``convolution'' on functionals:
\begin{equation}\label{E:convol}
\delta_{\Gamma_1}*\delta_{\Gamma_2}=\sum_\Gamma c_\Gamma\delta_\Gamma, \quad 
\Gamma_1\circ\Gamma_2=\sum_\Gamma c_\Gamma\ \Gamma.
\end{equation}
It will be shown that $*$ it is the convolution corresponding to the coproduct 
$\Delta_b$ introduced initially in \cite{pqft},
which is therefore dual to $*$.
\begin{defin}\label{D:coprod}
The {\em reduced graph coproduct} is:
\begin{equation}\label{E:coprod}
\Delta_b\Gamma=\sum_{\gamma\overset{j}{\hookrightarrow}\Gamma,\ \gamma\ \cap\ \partial\Gamma\neq\emptyset}
(-1)^{|j|}(\Gamma/\gamma)\otimes\gamma,
\end{equation}
where the sum corresponds to connected ``normal'' non-trivial subgraphs 
with a ``connected'' intersection with the boundary of $\Gamma$
($\partial\gamma$ consists of consecutive points of $\partial\Gamma$).
The sign is determined by $|j|$,
which denotes the number of boundary points to the left of $Im(j)$.
The corresponding {\em coproduct} will be denoted by:
$$\Delta=Id\otimes 1+1\otimes Id+\Delta_b.$$
\end{defin}
The ``normal'' subgraphs are those for which 
the quotient graphs belongs to the same class of graphs.
As noted before,
one may choose to further restrict the class $G$,
by discarding the graphs in the kernel of $\C{U}$
or with zero contribution in the definition of the star-product.

If $\Gamma\in G_{n,3}^{-1}$,
with $\partial\Gamma=\{1',2',3'\}$ (our case of interest), 
then the terms of $\Delta_b\Gamma$ belong to 
$G_{k,2}^0\otimes G_{l,2}^0$, where $k+l=n$.
The negative sign corresponds to collapsing the two boundary points
from the right $\{2',3'\}$,
while the positive terms correspond to subgraphs $\Gamma_2$
containing the boundary points $\{1',2'\}$.

Again a few examples of reduced coproducts are in order, 
dual to the above graph compositions:
\begin{equation}\label{E:co1}
\Delta_b(t_2^L)=b_2^L\otimes b_0-b_1\otimes b_1,\quad 
\Delta_b(t_2^R)=-b_2^R\otimes b_0+b_1\otimes b_1,
\end{equation}
\begin{equation}\label{E:co2}
\Delta_b(c_2^L)=b_1\otimes b_1-b_1^2\otimes b_0,\quad 
\Delta_b(c_2^R)=-b_1\otimes b_1+b_1^2\otimes b_0.
\end{equation}
The relation between the constructive definition \ref{D:star1}
and the concept of ``graph extension'' and the above coproduct,
is the content of the following lemma.
\begin{lem}
(i) For any graphs $\Gamma_1\in G_{k,2}^0, \Gamma_2\in G_{l,2}^0$:
$$\Gamma_1\circ\Gamma_2=
\sum_{\Gamma_2\overset{j}{\hookrightarrow}\Gamma\to\Gamma_1}
(-1)^{|j|}\Gamma,$$
where the sum is over the graphs $\Gamma\in G_{k+l,3}^{-1}$ such that
$\Gamma_2$ is a connected subgraph whose intersection with the boundary of $\Gamma$
is not empty and connected.

(ii) The operation $*$ \ref{E:convol}, dual to the pre-Lie operation $\circ$,
is the convolution corresponding to the coproduct \ref{E:coprod}:
\begin{equation}\label{E:pldual}
<\Gamma_1\circ\Gamma_2,\Gamma>=<\Gamma_1\otimes\Gamma_2,\Delta_b\Gamma>,
\end{equation}
\end{lem}
\begin{pf}
The coefficient of $\Gamma\in G_{n,3}^{-1}$ in 
$\Gamma_1\circ\Gamma_2$ is the algebraic sum of {\em signed extensions} $\pm\Gamma$
of $\Gamma_1$ by $\Gamma_2$,
such that $\Gamma_2$ is a subgraph of $\Gamma$ satisfying the restrictions imposed on
the coproduct \ref{E:coprod} and determining the sign $(-1)^{j}$,
while $\Gamma_1$ is the corresponding quotient.
This establishes the first claim (i).
Equivalently, (i) can be interpreted as stated in (ii)
by viewing the inner product as a duality and graphs as functionals,
concluding the proof.
\end{pf}
\begin{rem}
The requirement that any component of the subgraph should intersect the boundary
as a connected subset is necessary to ensure that the collapsing operation is well-defined,
and the order of the points on the boundary is preserved.
\end{rem}
\begin{defin}\label{D:bred}
The {\em $i^{th}$ boundary reduction} of $\Gamma\in G_{n,m}^l$ is:
$$\sigma_i(\Gamma)=(\Gamma/b_0^i)\otimes b_0,$$
where $b_0^i$ is the pair of boundary points $(i,i+1)$.
\end{defin}
For example (with $L/R$ notation for $1,2$):
\begin{equation}\label{E:sigma}
\sigma_L(\Gamma)=(\Gamma/b_0^L)\otimes b_0, \quad \sigma_R(\Gamma)=(\Gamma/b_0^R)\otimes b_0, 
\qquad \Gamma\in G_{n,3}^{-1},
\end{equation}
represent the ``top left degree'' terms of the coproduct,
and it is related with the ``merger operation'' from \cite{Kath}, p.17 (see Definition \ref{D:merger}).

The restriction of the reduced coproduct to prime graphs plays a special role
in defining the logarithm of the star-product
(sum over prime graphs).
\begin{defin}
The {\em Lie algebra of graphs} $\C{L}$ is the pre-Lie subalgebra of 
prime graphs (e.g. $b_1\in\C{L}$, but not $b_1^2$).
The {\em cobracket} of $\C{L}$ is:
\begin{equation}\label{E:cobra}
\Delta_p:\C{L}\to \C{L}\otimes \C{L}, \qquad 
\Delta_p\Gamma=\sum_{\C{L}\ni\gamma\overset{j}{\hookrightarrow}\Gamma,
\ \gamma\ \cap\ \partial\Gamma\neq\emptyset}
(-1)^{|j|}(\Gamma/\gamma)\otimes\gamma,
\end{equation}
where the sum is over prime graphs $\gamma\in\C{L}$.
\end{defin}
\begin{rem}\label{CK}
Note that a normal subgraph of a prime graph is prime
and the corresponding quotient is still a prime graph.
As a consequence $\Delta_p$ is the restriction of $\Delta_b$ and 
is therefore the dual of the pre-Lie operation $\circ$ on prime graphs $\C{L}$.

This case is similar to case of the Connes-Kreimer Hopf algebra \cite{CKha}, p.35-37)
where the dual of the pre-Lie operation is the part of the coproduct
corresponding to ``prime elements'' (simple cuts and trees),
so that the Hopf algebra is the UEA of the LA of ``prime elements''.
\end{rem}
As supporting evidence to the above relation we will reformulate the coproduct of $H_T$
in terms of subgraphs rather than using cuts.
\begin{lem}
The coproduct of $H_T$ has the following structure:
$$\Delta_{CK}(t)=\sum_{\gamma\subset t}(t/\gamma)\otimes\gamma, \quad t\in \C{T},$$
where the non-trivial subgraphs $\gamma$ are normal subgraphs 
whose components meet the boundary of the tree $t$ (its terminal points).
\end{lem}
\begin{pf}
The terms of the coproduct correspond to cuts which define
subgraphs (forests).
Excepting the empty set and $t$ itself,
the proper subgraphs are {\em normal} if the quotient is still a tree.
\end{pf}
The terms of the cobracket representing the reduction of the coproduct 
to the Lie algebra of prime graphs correspond to simple cuts.

Returning to the study of the pre-Lie operation $\circ$ itself,
the crucial property of the ``dual Feynman rule'' $\C{U}$
is the compatibility between the composition of (labeled) graphs and 
associated operators (Gerstenhaber composition).
\begin{lem}\label{L:pre-Lie}
$\C{U}$ is compatible with the pre-Lie compositions:
$$\C{U}_{\Gamma_1})(\xi_1)\circ \C{U}_{\Gamma_2}(\xi_2)=
\C{U}_{\Gamma_1\circ\Gamma_2}(\xi_1\wedge\xi_2).$$
\end{lem}
\begin{pf}
It is enough to note the corresponding relation at the level of individual
insertions of labeled graphs:
$$<\C{U}_{\Gamma_1}(\xi_1)\circ_i\C{U}_{\Gamma_2}(\xi_2), f>=
<\C{U}_{\Gamma_1\circ_i\Gamma_2}(\xi_1\wedge\xi_2),f>,\quad \Gamma_k\in G_{n_k,m_k},\ k=1,2,$$
where $f\in A^{m_1+m_2-1}$ and $\xi_i\in T(\wedge A)$ (see \cite{cfg} for details).
\end{pf}
We will need here (unlabeled graphs case) the above result for the special case
when the same fixed tensor is associated to all the internal vertices.
\begin{defin}\label{D:kr}
A {\em Kontsevich rule} with {\em propagator} $\alpha$
(a fixed tensor of rank $n$),
associates to an n-valent \underline{unlabeled} graph $\Gamma\in G_{n,m}$
the Hochschild cochain $\C{U}_\Gamma(\alpha^n)$,
also denoted $\C{U}(\Gamma)$.
\end{defin}
With this restriction and notation, 
the above Lemma specializes to the following.
\begin{lem}\label{L:kr}
A Kontsevich rule is a morphism of pre-Lie algebras:
$$\C{U}:(H,\circ)\to (H_{Hoch}(A,A),\circ), \quad \C{U}(\Gamma_1\circ\Gamma_2)=\C{U}(\Gamma_1)\circ\C{U}(\Gamma_2).$$
\end{lem}
From the above lemma,
the pullback of the associativity condition for the star-product follows.
\begin{prop}\label{P:cc1}
The star-product \ref{E:star-product} is associative $\star\circ\star=0$ iff:
\begin{equation}\label{E:assoc2}
W(\Delta_b(\Gamma))=0 \quad \Gamma\in G_{n,3}^{-1}.
\end{equation}
\end{prop}
\begin{pf}
Since $\C{U}$ is linear,
equation \ref{E:assoc1} reduces conform Lemma \ref{L:pre-Lie} 
to a constraint equation at the level of graphs:
$$0=\sum_{k+l=n}\sum_{\Gamma_1\in G_{k,2}, \Gamma_2\in G_{l,2}} 
W_{\Gamma_1} W_{\Gamma_2} 
\Gamma_1\circ\Gamma_2=
\sum_{\Gamma\in G_{n,3}^{-1}} c_\Gamma \Gamma,$$
with coefficients determined by $W$ via the coproduct:
\begin{alignat}{2}
c_\Gamma&=\sum_{k+l=n}\sum_{\Gamma_1\in G_{k,2}, \Gamma_2\in G_{l,2}} 
W_{\Gamma_1} W_{\Gamma_2} <\Gamma_1\circ\Gamma_2,\Gamma>\\
&\sum_{k+l=n}\sum_{\Gamma_1\in G_{k,2}, \Gamma_2\in G_{l,2}} 
W_{\Gamma_1\otimes \Gamma_2} <\Gamma_1\otimes\Gamma_2,\Delta_b(\Gamma)>\\
&=\cdot \circ (W\otimes W) (\Delta_b(\Gamma))=W(\Delta_b(\Gamma)),
\end{alignat}
where in the last equality the multiplicative
extension of $W$ to the tensor algebra of $H$ is still denoted $W$.
\end{pf}
Before investigating this ``cocycle condition'' \ref{E:assoc2}
in the section \ref{S:ca} (see also \cite{pqft}),
we will comment on some consequence relevant to the deformation theory
of associative algebras.

\section{The DGLA of graph deformations}\label{S:dgla}
The pre-Lie algebra of graphs $(H,\circ)$ has a distinguished element
$b_0$ which is ``associative'': $b_0\circ b_0=0$.
It therefore defines a differential $\delta=[b_0,\ ]$,
where $[\ ,\ ]$ is the corresponding graded Lie bracket.

The situation is completely analogous with the Hochschild
DGLA which controls the deformation theory
of associative algebras.
This is due to the fact that graph composition
is the pullback of Gerstenhaber composition
via the ``Kontsevich rule'' $\C{U}$ (a kind of a dual Feynman rule).
 
This approach follows the current trend of interpreting
the linear structures as representations of some discrete
geometric objects/categories,
e.g. algebras as being representations of operads 
(bialgebras and PROPs, cobordism categories and TQFTs etc.)

The essential point is that the Kontsevich rule
is a cohomological pairing,
associating to graphs and states on graphs polydifferential operators
(similar to Feynman rules and Feynman integrals):
$$\C{U}_\Gamma(\xi)=<\Gamma,\xi>, \quad
<\delta\Gamma,\xi>=d_{Hoch}<\Gamma,\xi>, \qquad \Gamma\in G,\ \xi\in \Lambda T.$$
\begin{defin}
The {\em graph pre-Lie algebra} is $(\C{L},\circ)$.
The associated {\em graph deformation DGLA} is $(\C{L},[\ ,\ ],d)$,
with differential $d=[b_0,\ ]$ determined by the {\em unit} $b_0$ of $\C{L}$.
\end{defin}
The graph deformation DGLA $(H,[,],\delta)$ 
controls the deformation of associative algebras.
The {\em Kontsevich representations} of $H$
(Kontsevich rules for a fixed 2-tensor $\alpha$, not necessarily Poisson: Definition \ref{D:kr}),
are the Hochschild DGLAs.
The graph deformation DGLA is the dual picture
to the Hopf algebra of graphs ($\circ$ dual to $\Delta_p$),
which is essentially the UEA of $\C{L}$ (see \cite{CKha}).

In view of the deformation theory,
a star-product corresponds to differentials $D=[Z, ]$:
$$D^2=0 \quad \Leftrightarrow\quad Z\circ Z=0,$$
which are perturbations of $\delta$.
Equivalently, if $Z=b_0+Z^+$, 
then $Z^+$ is a solution of the Maurer-Cartan (deformation) equation in $H$:
\begin{equation}\label{E:MC}
\delta Z^++\frac12 [Z^+,Z^+]=0.
\end{equation}
\begin{th}\label{T:deform}
Associative deformations are Kontsevich representations of deformations of the unit $b_0$ 
in the direction of the universal cocycle $b_1$ (Equation \ref{E:comp1}):
$$Z=b_0+b_1\epsilon+...,\qquad [b_0,b_0]=0, \quad \delta(b_1)=[b_0,b_1]=0,$$
satisfying the Maurer-Cartan equation \ref{E:MC}.
\end{th}
In view of the above cohomological pairing,
a star-product may be represented as:
$$\star=<Z,\ exp(\alpha)>.$$
Of course $[b_1,b_1]\ne 0$ ($b_1$ is not a solution of Maurer-Cartan equation), 
but it will be shown that there is no ``primary obstruction''
for the deformation equation:
$$b_1\circ b_1=\delta Z_2,$$
and that in the linear case,
with no loops (at tree-level),
there is a canonical ``semi-classical'' solution,
we call the {\em Moyal element}.
The coefficients of its logarithm $X=\ln Z$ 
(the {\em Hausdorff element} - expansion over prime graphs $\Gamma\in\C{L}$)
are expected to be the universal coefficients of 
the Hausdorff series of a Lie algebra.

Properties of the Bernoulli graphs,
graph composition and its dual coproduct
will be investigated next.

\section{The combinatorial approach}\label{S:ca}
Linearity of $\C{U}$ allows to pullback the associativity condition
of $\star$ at the level of graphs:
$$\star=\sum W_\Gamma\C{U}_\Gamma=\C{U}_{\sum W(\Gamma)\Gamma}=\C{U}(T),$$
where $Z$ is the dual of the character $W$:
$$Z=\sum_{n\ge 0}Z_n\epsilon^n, \qquad Z_n=\sum_{\Gamma\in G_{n,2}^0} W(\Gamma)\Gamma.$$
Then $\star\circ \star=0$ translates as $Z\circ Z=0$, which in terms of 
the coproduct is the equation from proposition \ref{P:cc1},
or in terms of the associated convolution: $W^2=0$.

We will start with an example specializing the above statements
to the case of the star-product modulo terms quadratic in $\epsilon$,
and rephrase the correction analysis from \cite{Kath}, \S 3.2,
using our tools introduced above.
In improvement also worth to be noted 
consists in using the Gerstenhaber bracket $[*,*]$,
rather then the individual terms, e.g. $*\circ_1 *$.

It will be also useful to associate a planar forest to graphs without loops,
by ``cutting'' the boundary points $L$ and $R$.
The corresponding operation,
as well as its linear extension to linear combinations of graphs, 
will be denoted as $C:G_{n,2}^l\to H$.
The (no loops) graph is prime iff the associated forest is a tree.

At this point we will introduce the {\em global super symmetry} on graphs:
\begin{equation}\label{E:symm}
T:H\to H, \quad T(\Gamma)=-\Gamma^t,
\end{equation}
where {\em transpose} graph $\Gamma^t$ is obtained by flipping the graph 
in the plane to revers the order of the boundary points.
In our special case $m=2$,
this amounts to the exchange of the ``left'' and ``right'' boundary points.
For example:
$$T(b_2^L)=-b_2^R, \ T(c_2^L)=-c_2^R.$$
The global symmetry is an involution $T^2=Id_H$
compatible with the coproduct and associated convolution $*$ 
(composition $\circ$).
\begin{lem}
(i) $T(\Gamma_1\circ\Gamma_2)=T(\Gamma_1)\circ T(\Gamma_2)$,
(ii) $T(\Gamma_1\Gamma_2)=-T(\Gamma_1)T(\Gamma_2)$,
(iii) $T(\Delta\Gamma)=\Delta(T\Gamma)$.
\end{lem}
\begin{pf}
(Direct check).
\end{pf}
For example (Equations \ref{E:co1}, \ref{E:co2}):
$$T(\Delta(t_2^L))=\Delta(t_2^R), \ T(\Delta(c_2^L))=\Delta(c_2^R).$$

%
%
\subsection{Correction analysis}\label{S:coa}
To exemplify the operations introduced so far
we rephrase the correction analysis from \cite{Kath}.

Consider the linear part of the star-product, i.e.:
$$Z_{<2}=Z_0+\epsilon Z_1,\quad
Z_0=w_0\ b_0, \ Z_1=w_1\ b_1.$$
Then modulo quadratic terms (with notation $Z^{\circ 2}=Z\circ Z$ etc.):
$$Z_{<2}^{\circ 2}\cong Z_0^{\circ 2}+[Z_0,Z_1]\epsilon=
w_0^2(b_0\circ b_0)+w_0 w_1(b_0\circ b_1+b_1\circ b_0)\epsilon=0$$
where we used the previous computations (Equations \ref{E:comp1}).

Therefore irrespective of the first two coefficients,
$Z^{\circ 2}$ starts with quadratic terms.
It is natural to normalize the series in order to correspond
to a deformation of multiplication in the direction of the Poisson bracket:
$w_0=w_1=1$.

%
%
\subsubsection{Constant Poisson structure}
As a simpler case to investigate,
consider first a constant Poisson structure (see \cite{Kath}, p.16).
Then $\C{U}$ vanishes when evaluated on $t_2^{L/R}$,
and (Equation \ref{E:comp1}):
$$Z_{<2}^{\circ 2}=\epsilon^2 w_1^2(c_2^L-c_2^R).$$
To ``compensate'' this term,
one has to add the graph $b_1^2$,
obtained ``merging'' two of the three boundary vertices of $c_2^{L/R}$.

This prompts for an operation ``inverse to cutting'' the boundary vertices,
reducing their number by one. 
%
\begin{defin}\label{D:merger}
The {\em merger operation} on graphs, $\sigma_m:G_{n,3}^{-1}\to G_{n,2}^0$,
is defined by:
$$\sigma_m(\Gamma)=\Gamma/{b_0^L}-\Gamma/b_0^R,$$
where $b_0^L=b_0 \cup \bullet$ ($b_0^R=\bullet\cup b_0$) 
consists from the two left (right) boundary points,
while $\bullet\in G_{0,1}^{1}$ is the unit for $\circ$.
\end{defin}
Note that the merger operation
is related with the boundary reduction operation (Definition \ref{D:bred}).

%
\begin{rem}\label{R:rec1}
This operation will be relevant to the recursive formula for the coefficients of $W$,
since it corresponds to the ``leading term'' in the coproduct:
$$\Delta_b\Gamma=\sum \gamma\otimes\gamma'+
(\Gamma/b_0^L)\otimes b_0-(\Gamma/b_0^R)\otimes b_0
=\sum\gamma\otimes\gamma'+\sigma(\Gamma)\otimes b_0,$$
where the degrees of $\gamma$ and $\gamma'$ are smaller then the degree of $\Gamma$.
The main point to keep in mind (to be established later) 
is that one of the last two terms belongs to the kernel of $\C{U}$,
allowing to solve for the other corresponding coefficient.
\end{rem}
%
\begin{rem}\label{R:CK}
Our graphs without circuits correspond to binary forests under the operation
of gluing all left and right terminals respectively:
$$\diagram
H_T \dto_{\sigma_T} & \\
G_{n,2}^0 & G_{n,3}^{-1}\lto^{\sigma}.
\enddiagram$$
The relation with the Connes-Kreimer Hopf algebra of rooted trees
(with the corresponding coefficients) will be investigated elsewhere.
\end{rem}
The correction term to be added to $Z_{<2}$ must be a multiple of $\sigma(c_2^L)=b_1^2$,
and corresponds to the coproduct:
$$\Delta_b(c_2^{L})=b_1\otimes b_1-b_1^2\otimes b_0.$$
Let $Z'=Z_{<2}+w_2b_1^2$.
Then the new quadratic terms produced are:
$$[b_0,b_1^2]=-(c_2^L-c_2^R)$$
and modulo cubic order terms:
$${Z'}^{\circ 2}\equiv (w_1^2b_1^2+w_0 w_2 [b_0,b_1^2])\epsilon^2=
\epsilon^2 (w_1^2-w_0\ w_2)(c_2^L-c_2^R).$$
With $w_2=1$,
$Z'$ gives an associative product up to cubic terms.

Note that the multiplicativity of $W$ requires $W(b_1^2)=W(b_1)^2$,
and the equation $W(\Delta_bc_2^L)=0$ is already satisfied
since $W(\Delta_b(b_1^L))=0$:
$$W(\Delta_b c_2^L)=W(\Delta_b(b_1^L))\ W(\sigma(b_1^M))=0.$$
The independent equations $W(\Delta_b\Gamma)=0$ seems to correspond to
$\Gamma\in G_{n,3}^{-1}$ a non-symmetric prime graph (see Lemma \ref{L:symm}),
but at this point it is not clear if there are any ``obstructions''.
\begin{rem}
As mentioned at the beginning, a symmetry factor of 2 was lost
in the above computation.
With a labeled $b_1^2$, $[b_0,b_1^2]$ would be twice $c_2^R-c_2^L$
(see Section \ref{S:symm}).
\end{rem}
In the constant coefficient case,
at the $n^{th}$ level stage the needed correction term is $b_1^n$ 
($\frac1{n!}b_1^n$ for labeled graphs),
yielding the exponential function in the Moyal product \cite{Kath}, p.18
(see also Section \ref{S:cps}).

\begin{rem}\label{R:kernel}
Note that for constant or linear Poisson coefficients,
one should work at the level of graphs modulo the kernel of $\C{U}$.

Also note that the correction terms appearing via merger of boundary points
form a subclass of graphs (posets),
corresponding to the class of binary trees.
More precisely,
they belong the class of graphs $\C{T}$
of graphs generated by forests of binary trees under
the merger and composition operations: $\sigma, \circ$ (see Remark \ref{R:CK}).
\end{rem}
%
%
%
\subsubsection{General Poisson structure}
In the general case, 
the coefficient of the quadratic (first non-zero) term of $Z\circ Z$ is:
\begin{equation}\label{E:quadr}
Z_1^{\circ 2}+[Z_0,Z_2]
=w_1^2\ b_1^{\circ 2}+ w_0w_2[b_0,b_1^2] 
+ w_0w_2^L[b_0,b_2^L]+w_0w_2^R[b_0,b_2^R],
\end{equation}
where:
$$Z_2=w_2^2b_1^2+w_2^L b_2^L+w_2^R b_2^R.$$
A direct computation gives:
$$[b_0,b_2^L]=t_2^L+c_2, \quad [b_0,b_2^R]=-(t_2^R+c_2).$$
To cancel the common symmetric term $c_2$,
one must choose symmetric coefficients: 
$$w_2'=w_2^L=w_2^R.$$
Then equation \ref{E:quadr} becomes:
$$(w_1^2-w_0w_2)(c_2^L-c_2^R)+w_1^2(t_2^R-t_2^L)+w_0w_2'(t_2^L-t_2^R).$$
Assuming $W$ multiplicative and normalized, and therefore $w_2=w_1^2=1$ and $w_0=1$,
the equation \ref{E:quadr} simplifies to:
$$(w_1^2-w_0w_2')(t_2^R-t_2^L).$$
It vanishes if $w_2'=1$,
and therefore $Z$ yields a star-product modulo cubic terms.

Note that the above correction did not use the Jacobi identity (see \cite{CI}, p.20):
$$Alt(U_{c_2}(\alpha^2))=0,$$
which at the level of graphs is:
\begin{equation}\label{E:jacobi}
b_1\circ b_1=c_2+c_2^L-c_2^R,
\end{equation}
as a direct computation shows.
Taking into account \ref{E:comp1} it is equivalent to:
\begin{equation}\label{jac2}
t_2^R-t_2^L=c_2.
\end{equation}
% ...
\begin{rem}
Note that the failure of $b_1$ to satisfy Maurer-Cartan equation:
\begin{equation}\label{E:MC1}
\delta(b_1)+1/2[b_1,b_1]=b_1\circ b_1,
\end{equation}
is mapped by $\C{U}$ to the associator of $\alpha$.
The Jacobi identity in turn can be written using the ``graph associator''
under Gerstenhaber composition without incorporating the Leibniz rule:
$$\alpha_g=c_2.$$
\end{rem}
As noted in Remark \ref{R:kernel},
one may look for coefficients $w$ such that $Z\circ Z\equiv 0$ modulo
the kernel of $\C{U}$,
which contains the above relation.
In what follows we will look for solutions without the above relation,
i.e. look for generic associative deformations of associative operations in 
an arbitrary direction.

%
%
\subsection{Symmetry factors}\label{S:symm}
To restore the symmetry factors,
consider the map $F:H_L\to H$ forgetting the labels of graphs
of the Hopf algebra of labeled graphs.
In order to match graph composition $\circ_L$ of $H_L$
with graph composition $\circ$ of $H$,
the following scaling must be applied.
\begin{lem}
If $F(\Gamma,l)=\Gamma|Aut(\Gamma)|/deg(\gamma)!$ then
$F$ is a morphism of Hopf algebras.
\end{lem}
For example:
$$F(\tilde{b}_1^2)=b_1, \quad F(\tilde{c}_2^L)=1/2c_2^L,$$
where $\tilde{\Gamma}$ denotes a labeled version of the graph $\Gamma$.
Then the composition of labeled graphs exhibits some symmetry factors.
For instance:
$$ [\tilde{b}_0,\tilde{b}_1^2]=
\tilde{b}_0\circ_L \tilde{b}_1^2+\tilde{b}_1^2\circ_L b_0=
2(\tilde{c}_2^R-\tilde{c}_2^L).$$

The correct coefficients for a star-product may be obtained as follows:
$$Z=\sum \frac{W(\Gamma)}{Aut(\Gamma)}\Gamma,$$
where $W$ is a solution of $W*W=0$ in $H^*$.

For example, modulo cubic terms, and omitting $\tilde{}$:
$$Z=b_0+b_1\epsilon+(b_1^2/2!+b_2^L+b_2^R)\epsilon^2+...,$$
$$b_1\circ b_1+\frac{1}{2} [b_0,b_1]=t_2^R-t_2^L,$$
$$[b_0,b_2^L]+[b_0,b_2^R]=t_2^L-t_2^R,$$
$$Z\circ Z\equiv(b_1\circ b_1+\frac12[b_0,b_1^2]+[b_0,b_2^L]+[b_0,b_2^R])\epsilon^2=0.$$ 
Therefore the above $Z$ defines an associative star-product
up to cubic terms.

The above correction analysis has similarities with 
the BPHZ renormalization procedure,
both relying at some point on the Hopf algebra of trees
as a device encoding the combinatorics 
of graph composition.

%
%
\subsection{Is there a recursive formula?}\label{S:rf}
The first few coefficients of $W$ 
may be found solving the cocycle condition $W(\Delta_b\Gamma)=0$.

Since the Bernoulli graphs $b_n$ are primitive elements,
$W(\Delta(b_i))=0$, and there are no constraints for $n=0,1$.
We need $W(b_0)=1$ and the normalization 
$W(b_1)=1$ will be assumed for simplicity.

Writing explicitly the reduced coproduct for prime graphs 
with $n=2$ yields the previous values determined using 
``correction analysis'':
$$W(\Delta_b(t_2^L))=0, \quad W(-b_1\otimes b_1+b_2^L\otimes b_0)=0,
 \quad => \quad W(b_2^L)=1,$$
$$W(\Delta_b(c_2^L))=0, \quad W(b_1\otimes b_1-b_1^2\otimes b_0)=0,
 \quad => \quad W(b_1^2)=1.$$
%

To understand the structure of the coproduct,
introduce the following.
\begin{defin}
For any $\Gamma\in G_{n,2}^0$,
let $\phi_L(\Gamma)$ ($\phi_R(\Gamma)$) 
denote the number of paths in $\Gamma$ 
ending at the left (right) boundary point ``L'' (respectively ``R'')
of $\Gamma$.

The {\em hight of $\Gamma$} is $\phi(\Gamma)=\phi_R(\Gamma)-\phi_L(\Gamma)$.
\end{defin}
For example $\phi_L(b_2^L)=2,\ \phi_R(b_2^L)=1$, so $\phi(b_2^L)=-1$.
In general $\phi(b_n^L)=1-n$.
\begin{lem}
For all $\Gamma\in G_{n,3}^{-1}$:
$$\Delta_b(\Gamma)=
\sum\gamma_1\otimes\gamma_2+\
(\Gamma/b_0^L)\otimes b_0-(\Gamma/b_0^R)\otimes b_0,$$
where $deg(\gamma_i)<deg(\Gamma)$ and $b_0^{L/R}$ are the subgraphs 
consisting of the two left/right boundary points of $\Gamma$.
Moreover, if both terms are present, then:
\begin{equation}\label{E:h}
\phi_R(\Gamma/b_0^L)<\phi_R(\Gamma/b_0^R).
\end{equation}
\end{lem}
\begin{pf}
Note that there are only two subgraph $\gamma\subset\Gamma$
without interior points: $b_0^{L/R}$.
And any other subgraph will have a positive degree,
decreasing accordingly the degree of the quotient.

The last part follows from the definition of $\phi_R$,
since if the terms are not equal,
there will be at least one path in $\Gamma$ ending at the
middle point of its boundary.
\end{pf}
Note that the quotients $\Gamma/(b_0^{L/R})$ are just the graphs
obtained by ``merging'' boundary points in the correction analysis section.
\begin{lem}\label{L:symm}
Let $\Gamma\in G_{n,3}^{-1}$.

(i) If $\Gamma$ is symmetrical ($\Gamma^t=\Gamma$),
then $\Delta_b\Gamma$ is anti symmetrical.
Consequently, if $W$ is equivariant, then:
\begin{equation}\label{E:W}
W(\Delta_b\Gamma)=0.
\end{equation}

(ii) If $\Gamma$ has no ``middle legs'' (arrows lending on the second boundary point),
then $\Delta_b\Gamma=0$.
Otherwise \ref{E:h} holds and \ref{E:W} may be solved for the coefficient $W$ 
of $\sigma_R(\Gamma)$.
\end{lem}
\begin{pf}
(i) The subgraphs $\gamma\subset\Gamma$ containing boundary points 
cannot be symmetrical,
as have to contain either ``L'' or ``R'',
but not both.
Therefore the terms of the coproduct come in mirror pairs
with opposite signs.
When applying a symmetrical $W$, the result is zero.

(ii) Follows by inspection (see \ref{E:sigma}).
\end{pf}
For example $\Delta_b(c_2)=b_2^L\otimes b_0-b_2^R\otimes b_0$,
and $W(\Delta_b(c_2))=0$ will not impose a constraint on 
the coefficients $W$.
\begin{prop}\label{P:bern}
The coefficients of a star-product at tree-level 
are determined by the coefficients of the Bernoulli graphs.
\end{prop}
\begin{pf}
By induction on the number of internal vertices
and ``hight function'' $\phi_R$
the coefficients $W$ are determined by the values
on the graphs which minimize the height function:
the left Bernoulli graphs.
\end{pf}
To investigate existence and uniqueness of a ``semi-classical'' solution
(graphs without loops / ``tree level''),
we will consider the relation with the antipode,
aiming for an explicit recursion relation.

%
%
%
\subsection{Relations with the antipode}
The associativity condition may be interpreted 
in terms of the antipode,
establishing a relation (to be explored) 
with the renormalization procedure using the ``counter terms'' method 
(BPH and Zimmermann's formulas, as well as Dyson-Salam procedure \cite{FG}).

Let $(H^*,*)$ be the group of characters of $H$,
equivariant with respect to the global symmetry $T$.
The inverse of a character $w$ is $S^*(w)=w\circ S$,
where $S$ is the antipode of $H$.

A few values of the antipode can be computed readily 
using the {\em recursive formula}:
$$S(\Gamma)=-\Gamma-\sum S(\Gamma/\gamma)\gamma,$$
where the terms in the above sum corresponds to the reduced
coproduct.
For the primitive elements $b_n^{L/R}$, $S(b_n^{L/R})=-b_n^{L/R}$.
Moreover, for prime elements with $n=2$ internal vertices we have:
$$S(t_2^L)=-t_2^L-b_1\otimes b_1+b_2^L\otimes b_0,\quad
S(c_2^L)=-c_2^L+b_1\otimes b_1-b_1^2\otimes b_0.$$
The right handed versions of the above formulas may be obtained
using the global symmetry: 
$$S(T(\Gamma))=T(S(\Gamma)).$$
%
There is an {\em explicit formula} for the antipode,
which is obtained by expressing the antipode as a geometric series \cite{FG}, p.3:
$$S=(id)^{-1}=\sum_{n\ge 0}x^{*n}, \quad x=\B{1}-id,\ \B{1}=u\circ\eta,$$
where $u$ is the unit of $H$ and $\eta$ the counit.
The unit in $End(H)$ is $\B{1}$ and $x^0=\B{1}$.

If $W$ is extended to a character $W_c=\eta+W$ on $H$ ($\eta$ counit),
then the associativity condition for the star-product \ref{E:assoc2}
corresponds to:
$$W(S(\Gamma))=-W(\Gamma), \quad i.e.\quad W_c^{-1}=\bar{W}_c,$$
for instance by using the recursive relation for $S$,
where $\bar{W}_c(\Gamma)=\eta-W_c(\Gamma)$ is the ``conjugate character''.

%
\section{Special cases revisited}\label{S:scr}
The constant and linear Poisson structures will be revisited.

%
\subsection{Constant Poisson structures}\label{S:cps}
In the constant Poisson structure case the ``constraint equation'' 
is identically satisfied and the well-known Moyal formula is recovered.

Recall that all internal vertices of graphs in $G_{n,2/3}^{0/1}$
have only two descendents ($\alpha$ is a degree two tensor),
but may have several ``parents'' ($\partial_i$s acting on the corresponding
coefficient of $\alpha$).

For a constant Poisson structure $\alpha=\sum c^{ij}\partial_i\wedge \partial_j,$
the non-trivial graphs $\Gamma\in \C{G}$ (in $G_{n,2}^0$, but not in the kernel of $\C{U}$) 
have bi-valent internal vertices with only outgoing edges,
therefore ending on boundary vertices:
$$\C{G}=\{\Gamma_n\}_{n\ge 0}, \ \Gamma_n=b_1^n, \ n\ge 0.$$
They are generated by the only non-trivial Bernoulli graph $b_1$.
The prime graphs of $G_{n,3}^{-1}$ are $b_1^{L/M/R}$, with coproducts:
$$\Delta_b(b_1^L)=b_0\otimes b_1-b_1\otimes b_0=-\Delta_b(b_1^R),
\quad \Delta_b(b_1^M)=0.$$
Then equation $W(\Delta_b\Gamma)=0$ is identically satisfied on ,
and $W\equiv1$ is determined by our initial conditions $w_0=w_1=1$.
With explicit symmetry factors the Moyal element is: 
$$Z=\sum_{\Gamma\in \C{G}}W(\Gamma)\Gamma=\sum_{n\ge 0}b_1^n/n!\hbar^n=exp(b_1\hbar).$$
The corresponding Hausdorff element (its logarithm) is $X=\hbar b_1$.
Representing $Z$ yields Moyal's formula: $\star=\C{U}(Z)$.

This is the commutative case corresponding to the abelian Lie algebra
with one generator $b_1$.
The ``general case'' from the point of view of Lie bialgebras \cite{Resh}
and corresponding to linear Poisson structures,
will be studied next.

%
\subsection{Linear Poisson structures}\label{S:lps}
In the linear Poisson structure case (e.g. dual of a Lie algebra \cite{Kath}, p.5)
the equation can be solved recursively,
yielding a canonical semi-classical (tree-level) star-product.
The coefficients are essentially those of the Hausdorff series
(CBH-quantization \cite{Kath}).

For a linear Poisson structure:
$$\alpha=\sum c^{ij}_k X^k \partial_i\wedge \partial_j,$$
the non-trivial graphs have only tri-valent internal vertices,
corresponding to trees.
\begin{lem}
The number of edges of $\Gamma\in G_{n,2/3}^{0/1}$ is $2n$.

If $\Gamma\in G_{n,2}^0$ is such that $\C{U}_\Gamma(\alpha^n)\ne 0$,
then each internal vertex has only one parent.
Equivalently, if $\Gamma$ has at least one n-valent internal vertex with $n>3$,
then $\C{U}_\Gamma(\alpha^n)=0$.
\end{lem}
\begin{pf}
Internal vertices with more then one parent will contribute a zero factor
to the corresponding product in any term of the state-sum.
\end{pf}
The correspondence between trees and the relevant class of graphs $\C{G}$
not belonging to the kernel of $\C{U}$:
\begin{equation}\label{E:corresp}
t:\C{T}\to \C{G}, \quad s:\C{G}\to \C{T}, \quad s=t^{-1}
\end{equation}
is defined by collecting together the left terminal points to the boundary point ``L''
and all right terminal points to the boundary point ``R''.
It extends to a linear map between the corresponding Hopf algebras $H_T$ and $H_G$.

It has an inverse $s$ (spanning tree/forest),
obtained by ``cutting'' the boundary points ``L'' and ``R''
(letting lose the terminal edges).

Now the (left) Bernoulli graphs,
corresponding to left parenthesized expressions / trees,
minimize the height function. 
Moreover the star-product is determined by the value of $W$ 
on Bernoulli graphs.
These in turn are determined recursively by noting that
for each $b_n$ there is a graph $\Gamma_n\in G_{n,3}^{-1}$ such that
$$\Delta_b\Gamma_n=b_n^L\otimes b_0-b_{n-1}\otimes b_1.$$
For example $\Gamma_3$ is the following graph:
$$\diagram
\bullet\dto\rto & \bullet \dlto \rto & \bullet \dlto \dto \\
\bullet & \bullet & \bullet.
\enddiagram$$
As a consequence the values of $W$ are determined inductively.
Since only tree-like graphs have been considered,
the corresponding star-product may be thought off as
a semi-classical solution, 
i.e. not including quantum corrections due to loops/circuits.
This star-product corresponds to the Hausdorff series,
except our coefficients neglect the symmetry factors,
yielding $W(b_n)=1$ (!?).
Its logarithm (a sum over prime graphs),
is expected to correspond to the Hausdorff series.
\begin{th}
A linear Poisson structure
has a unique semi-classical solution of the
initial value problem for the star-product at tree-level
$$\star_{|\epsilon=0}=m, \quad \frac{d\star}{d\epsilon}_{|\epsilon=0}=\{,\}.$$
The coefficients of the associated Moyal element 
are determined inductively from the Equation \ref{E:assoc2}.
\end{th}
It is expected that the antipode may be used to derive 
further formulas for the coefficients.

\section{Conclusions}\label{S:concl}
Deformation quantization based on star-products
require finding coefficients which are cocycles in some sense.
To quantize any Poisson manifold,
Kontsevich used a simple string theory 
(non-linear sigma model on the Poincare disk)
to write explicit formulas for the (finite) Feynman integrals / amplitudes.
Then the resulting cocycle was used in the general case.

In the present article a combinatorial approach is sketched,
leading to a self contained mathematical solution.

It was shown that the cocycles reside 
in the deformation DGLA of graphs controlling the deformation theory
of associative algebras.
Star-products appear as ``Kontsevich representations'' elements $Z$
such that $Z\circ Z=0$.
When restricted to graphs without circuits
a canonical solution of the deformation equation exists (\cite{Kon1}): 
the ``Moyal element''.
In the special case of Poisson structures with constant
coefficients, its representation yields Moyal's formula.
For linear coefficients, it is expected to correspond to 
the Hausdorff series.

The Hopf algebra techniques have proven as a powerful 
conceptual interface to the combinatorial practices
of renormalization \cite{Kre1,CKren1}.
A direct connection with renormalization
is likely to be established via the Hopf algebra of graphs (Connes-Kreimer, Loday-Ronco etc.).
For a linear Poisson structure,
the class of non-trivial (binary) prime graphs $\C{G}$ (not in the kernel of the ``Feynman rule'' $\C{U}$)
is in 1:1 correspondence with binary trees (see Equation \ref{E:corresp}).
Note that in this hierarchic case (no multiple parents), 
there is no need for considering labeled trees,
as in the renormalization of Feynman graphs.
A correspondence between the corresponding Hopf algebras is expected to originate
from a correspondence at the level of pre-Lie algebras.

Further considerentions regarding the relations with renormalization will be addressed elsewhere.

\section{Appendix}\label{S:ex}
The first few compositions and coproducts are given below.
The right-handed version may be obtained by symmetry,
applying the Hopf involution $T\Gamma=-\Gamma^t$.

Recall that for symmetric graphs ($\Gamma^t=\Gamma$) $\Delta\Gamma$ is antisymmetric,
and therefore $W(\Delta\Gamma)=0$ for any equivariant $W$.

For products of elements, 
the corresponding values are determined by the values of prime elements,
using the properties of $\Delta_b$.

The ``modifiers'' $L/R$ are used to denote a shift: $\Gamma_L=\Gamma\bullet$.
For example:
$$b_{1L}^2=b_1^2\bullet, \quad, b_{1R}^2=\bullet b_1^2.$$

\vspace{.2in}
{\bf Compositions of some elements from $G_{n,2}^0$}
Prime elements:

\noindent Degree 0:
$$b_0\circ b_0=0.$$
Degree 1:
$$b_0\circ b_1=b_1^L-b_1^R,\quad
b_1\circ b_0=b_1^R-b_1^L, \qquad [b_0,b_1]=0.$$
Degree 2:
$$b_0\circ b_2^L=b_2^L\ \bullet-\bullet\ b_2^L, 
\quad b_2^L\circ b_0=\bullet\ b_2^L-b_2^L\ \bullet+t_2^L+c_2,
\qquad [b_0,b_2^L]=t_2^L+c_2.$$
$$b_1\circ b_1=t_2^R-t_2^L+c_2^L-c_2^R.$$
Non-prime elements:
$$b_0\circ b_1^2=b_1^2 \bullet-\bullet b_1^2, \quad
b_1^2\circ b_0=\bullet b_1^2-b_1^2 \bullet+c_2^R-c_2^L,
\qquad [b_0,b_1^2]=c_2^R-c_2^L.$$

\vspace{.2in}
{\bf Coproducts of some prime elements from $G_{n,3}^{-1}$}
Prime elements:

\noindent Degree 0:
$$\Delta_b (\bullet\bullet\bullet)=0.$$
Degree 1:
$$\Delta_b(b_1^L)=b_0\otimes b_1-b_1\otimes b_0,\quad \Delta_b(b_1^M)=0.$$
Degree 2:
$$\Delta_b(t_2^L)=b_2^L\otimes b_0-b_1\otimes b_1,\quad 
\Delta_b(c_2)=b_2^L\otimes b_0-b_2^R\otimes b_0.$$
Non-prime elements:
$$\Delta_b(c_2^L)=b_1\otimes b_1-b_1^2\otimes b_0,
\qquad \Delta(b_{1L}^2)=b_0\otimes b_1^2-b_1^2\otimes b_0, \qquad \Delta(b_{1M}^2)=0.$$
%



\end{document}